\documentclass[11pt]{amsart}
\usepackage[all]{xy}

\usepackage[utf8]{inputenc}		% font and encoding
\usepackage[T1]{fontenc}
\usepackage[english]{babel}

\usepackage{amsfonts}			% ams-packages
\usepackage{amsmath}
\usepackage{amssymb}
\usepackage{amsthm}
\usepackage{mathtools}
\makeatletter					% Add the 2020 MSC, which is not yet in the
\@namedef{subjclassname@2020}{	% package by default. Remove this when it's added.
	\textup{2020} Mathematics Subject Classification}
\makeatother

\usepackage{graphicx}			% images and drawing
\usepackage{tikz}
	\usetikzlibrary{cd}

\usepackage{hyperref}			% hyperlinks
	\hypersetup{colorlinks=false, urlcolor=black, linkcolor=black}

\usepackage{enumitem}			% arbitrary list items

\usepackage{color}				% colored text for signaling changes

\newcommand{\Z}{\mathbb{Z}}						% Integers
\newcommand{\R}{\mathbb{R}}						% Reals
\renewcommand{\S}{\mathbb{S}}					% Sphere
\newcommand{\eps}{\varepsilon}					% Epsilon shortcut
\newcommand{\dd}								% Differential d
	{\mathop{}\!\mathrm{d}}						
\newcommand{\ddn}[1]							% Powers of a differential d
	{\mathop{}\!\mathrm{d^{#1}}}

\newcommand{\abs}[1]							% Absolute value
	{\left| #1 \right|}
\newcommand{\smallabs}[1]						% Small absolute value bars which won't scale to the argument.
	{\lvert #1 \rvert}	
\newcommand{\norm}[1]							% Norm 
	{\left\lVert #1 \right\rVert}	
\newcommand{\smallnorm}[1]						% Small norm bars which won't scale to the argument.
	{\lVert #1 \rVert}						
\newcommand{\ip}[2]								% Inner product
	{\left< #1 , #2 \right>}

\DeclareMathOperator{\dist}{dist}					% Support
\newcommand{\loc}{\mathrm{loc}}					% Local spaces

\def\Xint#1{\mathchoice
	{\XXint\displaystyle\textstyle{#1}}%
	{\XXint\textstyle\scriptstyle{#1}}%
	{\XXint\scriptstyle\scriptscriptstyle{#1}}%
	{\XXint\scriptscriptstyle\scriptscriptstyle{#1}}%
	\!\int}
\def\XXint#1#2#3{{\setbox0=\hbox{$#1{#2#3}{\int}$}
		\vcenter{\hbox{$#2#3$}}\kern-.5\wd0}}

\def\dashint{\Xint-}

\newtheorem{thm}{Theorem}[section]{\bf}{\it}

\newtheorem{innercustomthm}{Theorem}[section]{\bf}{\it}
\newenvironment{customthm}[1]
	{\renewcommand\theinnercustomthm{#1}\innercustomthm}
	{\endinnercustomthm}
	
{\bf}{\it}

\newtheorem{innercustomlemma}{Lemma}[section]{\bf}{\it}
\newenvironment{customlemma}[1]
	{\renewcommand\theinnercustomlemma{#1}\innercustomlemma}
	{\endinnercustomlemma}

\theoremstyle{definition}

\theoremstyle{remark}

\newtheorem{innercustomrem}{Remark}[thm]
\newenvironment{customrem}[1]
{\renewcommand\theinnercustomrem{#1}\innercustomrem}
{\endinnercustomrem}

\numberwithin{equation}{section}

\begin{document}
	
\title{Corrigendum to ``On the heterogeneous distortion inequality''}
\author{Ilmari Kangasniemi}
\address{Department of Mathematical Sciences, University of Cincinnati, P.O.\ Box 210025, Cincinnati, OH 45221, USA}
\email{kangaski@ucmail.uc.edu}

\author{Jani Onninen}
\address{Department of Mathematics, Syracuse University, Syracuse,
NY 13244, USA and  Department of Mathematics and Statistics, P.O.Box 35 (MaD) FI-40014 University of Jyv\"askyl\"a, Finland}
\email{jkonnine@syr.edu}
	
\thanks{The authors I.\ Kangasniemi and J.\ Onninen are currently supported by the National Science Foundation grants DMS-2247469 and DMS-2154943, respectively.}
%\subjclass[2020]{Primary 30C65; Secondary 35B53, 35R45, 53C21}
%\keywords{Heterogeneous distortion inequality, Quasiregular mappings, Liouville theorem, H\"older continuity, Astala-Iwaniec-Martin question}

\begin{abstract}
We correct an error in [I. Kangasniemi, and J. Onninen, \emph{On the heterogeneous distortion inequality.}
Math. Ann. 384 (2022), no. 3-4, 1275–1308.]
\end{abstract}
	
\maketitle
	
When investigating the potential alternate applications of the methods in our paper \emph{On the heterogeneous distortion inequality}~\cite{Kangasniemi-Onninen_Heterogeneous}, we discovered that the presented proof of the main result, Theorem 1.3, has a critical flaw. This error occurs late in the paper, in the relatively technical proof of the lower integrability result shown in Lemma 7.2, and invalidates this lemma in its stated form. The only results in the paper affected by this error are Lemma 7.2 and Theorem 1.3.

We have been unable to reprove the original statement of Lemma 7.2. However, in this corrigendum, we present a fix that recovers the main result, Theorem 1.3, in its entirety. The fix is non-trivial, and took us numerous failed attempts to find.   

The error in the proof of Lemma 7.2 lies in the use of the Hardy-Littlewood maximal inequality. In particular, we apply it on a subset of $\R^n$. While there are maximal inequalities on more general domains, this requires that the definition of the maximal function is restricted to the domain; in our case, the definition of the maximal function extends past the domain, invalidating the inequality.

\begin{customrem}{1.1}
	We remark that a slightly weakened version of Theorem 1.3 follows from Theorem 1.2 of the original article with a far easier argument. In this version, instead of assuming that $\sigma \in L^{n+\eps}(\R^n) \cap L^{n-\eps}(\R^n)$, one would assume that $\sigma \in L^{n+\eps}(\R^n)$ and $(\sigma \circ \iota) \abs{D\iota} \in L^{n+\eps}(\R^n)$, where $\iota \colon \R^n \cup \{\infty\} \to \R^n \cup \{\infty\}$ is the conformal inversion across the unit sphere. This alternate version is weaker, as its assumptions imply $\sigma \in L^{n+\eps'}(\R^n) \cap L^{n-\eps'}(\R^n)$ for $\eps' \in (0,  n\eps/(n+2\eps))$. The proof essentially only amounts to pre-composing the map with $\iota$, proving that the resulting inverted map is in $W^{1,n}_\loc(\R^n)$ by using global $L^n$-integrability of the derivative and removability of points for $W^{1,n}$-spaces, proving a generalized distortion estimate for the inverted map, and then applying Theorem 1.2 to show that the inverted map cannot have a zero at the origin. We leave the details to the interested reader, as our following fix makes this weaker version unnecessary.
\end{customrem}

\section*{The corrected proof}

\setcounter{section}{7}

We follow the notational conventions of the original paper \cite{Kangasniemi-Onninen_Heterogeneous}. We start right before Lemma 7.2 in the original paper, where we have estabilished a continuous mapping $h \in W^{1,n}_\loc(\R^n, \R \times\S^{n-1})$ with $\abs{Dh} \in L^n(\R^n)$ and
\begin{equation}\label{eq:logdistineq}
	\abs{Dh(x)}^n \le K J_h(x) +\sigma^n(x) \qquad  K\in [1, \infty), \text{ a.e.\ } x\in \R^n \, .  
\end{equation}
We match our theorem numbering to the original article, and refer to results therein.

In the original flawed Lemma 7.2, the statement essentially said that if $\sigma \in L^n(\R^n) \cap L^{n-\eps}(\R^n)$, then $\abs{Dh} \in L^{n-\eps'}(\R^n)$ for some $\eps' \in (0, \eps')$. The general idea of this relied on adapting an argument by Faraco and Zhong~\cite{Faraco-Zhong_Caccioppoli} which is used to prove a Caccioppoli inequality below the natural exponent for mappings of bounded distortion. The difference in our case was that we aimed to perform the argument globally, leveraging the fact that $\abs{Dh} \in L^n(\R^n)$.

A na\"ive repetition of the argument of Faraco and Zhong in our setting is close to yielding an estimate of the desired form, but fails due to a potential $\infty - \infty$ -cancellation in one of the steps. The original argument was an attempt at a cut-off procedure that would make the relevant terms finite in the part of the proof where this cancellation occurs. However, we are currently unable to find any such cut-off methods that would fix the proof, largely since the Hardy-Littlewood maximal estimate that is used to undo an extra maximal function in the canceled term is global. 

Our solution is to instead perform a similar argument on a logarithmic scale of lower integrability. In contrast to the $L^p$-case, a difference occurs in the logarithmic case where the term to be absorbed has a logarithm less of lower integrability. This difference provides enough leverage to avoid the aforementioned issues in the cancellation step, by relying on a trick essentially based on interpolation of exponents. 

In particular, the technical estimate we show that replaces the previous Lemma 7.2 is as follows.

\begin{customlemma}{7.2}[revised]\label{prop:log_integrability}
	Suppose that a mapping $h \colon \R^n \to \R \times \S^{n-1}$ is continuous and non-constant, and that $h \in W^{1,n}_\loc(\R^n, \R \times \S^{n-1})$ with $\abs{Dh} \in L^n(\R^n)$.  If $h$ satisfies the distortion inequality \eqref{eq:logdistineq} with $\sigma \in L^{n-\eps}(\R^n) \cap L^{n}(\R^n)$ for some $\eps \in (0, n-1)$, then we have
	\[
		\int_{\R^n} \abs{Dh}^n \log^{n}\left( 1 + \frac{1}{M(\abs{Dh})} \right)
		\leq C(n, K, \eps) \int_{\R^n} (\sigma^n + \sigma^{n-\eps}) < \infty.
	\]
\end{customlemma}

We remark that the left hand side of the above estimate is indeed below the natural exponent on the log-scale, since $\log(1 + t^{-1}) \to \infty$ when $t \to 0$. The estimate is however weaker that the full $L^n \log^n(1 + L^{-1})$-estimate due to the maximal function in the logarithm, but it is regardless sufficient for us. 

For the proof of the revised Lemma \ref{prop:log_integrability}, given a positive integer $m \in \Z_{>0}$, we define the function $\Phi_m \colon (0, \infty) \to (0, \infty)$ by
\begin{align*}
	\Phi_m(r) 
	&= - \frac{d}{dr} \log^{m} (1 + r^{-1})\\
	&= \frac{m}{r(1+r)} \log^{m-1}(1+r^{-1}).
\end{align*}
Our key reason for using this function is that, since $\lim_{t \to \infty} \log^{m} (1 + t^{-1}) = \log^m(1) = 0$, we have
\begin{equation}\label{eq:Phi_integral_identity}
	\int_r^\infty \Phi_m(t) \dd t = \log^{m} (1 + r^{-1}).
\end{equation}
The other basic property of $\Phi_m$ that we need is the following integral estimate.

\begin{customlemma}{7.3}\label{lem:Phi_integral_lemma}
	For every positive integer $m \in \Z_{>0}$, we have for all $r > 0$ the estimate
	\[
		\int_0^r t \Phi_m(t) \dd t \leq C_m r \left( 1 + \log^{m - 1}(1 + r^{-1})\right).
	\]
\end{customlemma}
\begin{proof}
	Consider first the case $m = 1$. In this case, we can directly compute the integral and use the elementary estimate $\log(1+r) \leq r$ to obtain
	\[
			\int_0^r t \Phi_m(t) \dd t = \int_0^r \frac{dt}{1+t} = \log(1+r) \leq r = \frac{r}{2} \left( 1 + \log^0(1 + r^{-1}) \right).
	\]
	With the case $m = 1$ done, we then proceed by induction, supposing that $m \geq 2$ with the claim proven for $m - 1$. We first observe that for all $t > 0$, we have
	\[
		t\Phi_m(t) = \frac{m}{1+t} \log^{m - 1} (1 + t^{-1}) \leq m \log^{m - 1}(1 + t^{-1}).
	\]
	We then integrate by parts, differentiating $\log^{m-1}(1+t^{-1})$ to get $-\Phi_{m - 1}(t)$ and integrating $1$ to get $t$. We get the estimate
	\begin{multline*}
		\int_0^r t \Phi_m(t) \dd t
		\leq m \int_0^r \log^{m-1}(1+t^{-1}) \dd t\\
		= m \bigl(r \log^{m - 1} (1 + r^{-1}) - \lim_{s \to 0^{+}} s \log^{m - 1} (1 + s^{-1})\bigr) + m \int_0^r t \Phi_{m-1}(t) \dd t.
	\end{multline*}
	The limit term in the above upper bound vanishes, though it could regardless be ignored due to it clearly being negative. By using the induction assumption, we hence have
	\[
		\int_0^r t \Phi_m(t) \dd t \leq m r \log^{m - 1} (1 + r^{-1}) + m C_{m - 1} r \left( 1 + \log^{m - 2}(1 + r^{-1})\right).
	\]
	Since $m - 2 \geq 0$, we may then use an interpolation estimate to obtain that $\log^{m - 2}(1 + r^{-1}) \leq \max(\log^0 (1 + r^{-1}), \log^{m-1}(1 + r^{-1})) \leq 1 + \log^{m-1}(1+r^{-1})$, and the claim follows.
\end{proof}
\begin{proof}[Proof of the revised Lemma \ref{prop:log_integrability}]
	Note that we may in fact assume that $\sigma \leq (K+1)^{1/n} \abs{Dh}$. Indeed, since $\abs{J_h} \leq \abs{Dh}^n$ we always have for free the estimate $\abs{Dh}^n \leq KJ_h + (K+1)\abs{Dh}^n$. Thus, if $h$ satisfies \eqref{eq:logdistineq} with $\sigma = \sigma_0$, then $h$ also satisfies \eqref{eq:logdistineq} with $\sigma = \min(\sigma_0, (K+1)^{1/n} \abs{Dh})$.
	
	Somewhat similarly to the original Lemma 7.2, we define for every $\lambda > 0$ the set
	\[
		F_\lambda = \{x \in \R^n : M(\abs{Dh}) \leq \lambda\}.
	\]
	If $x, y \in F_\lambda$, we then again have by a pointwise Sobolev estimate that
	\begin{multline*}
		\abs{h_{\R}(x) - h_{\R}(y)} \leq (C_n/2) \abs{x - y} (M(\abs{\nabla h_{\R}})(x) + M(\abs{\nabla h_{\R}})(y))\\
		\leq (C_n/2) \abs{x - y} (M(\abs{Dh})(x) + M(\abs{Dh})(y)) \leq C_n \lambda \abs{x - y},
	\end{multline*}
	proving that $h_{\R}$ is $C_n \lambda$-Lipschitz in $F_\lambda$. By using McShane extension, we can then again find a $C_n \lambda$-Lipschitz map $h_{\R, \lambda} \colon \R^n \to \R$ such that $h_{\R, \lambda} \vert_{F_\lambda} = h_{\R} \vert_{F_\lambda}$. We again denote $h_\lambda = (h_{\R, \lambda}, h_{\S^{n-1}})$. 
	
	We continue as in the original argument, proving that 
	\[
		\abs{D h_\lambda} \leq (1+ C_n) M(\abs{Dh}).
	\]
	Indeed, in $F_\lambda$ we have $\abs{D h_\lambda} = \abs{Dh}$, and in $\R^n \setminus F_\lambda$ we have $\abs{D h_\lambda} \leq \abs{Dh} + C_n \lambda \leq (1+C_n) M(\abs{Dh})$. By the Hardy-Littlewood maximal inequality, we have $M(\abs{Dh}) \in L^n(\R^n)$, and thus \cite[Lemma 2.4]{Kangasniemi-Onninen_Heterogeneous} of the original paper yields that
	\[
		\int_{\R^n} J_{h_\lambda} = 0
	\]
	for every $\lambda > 0$.
	
	In $\R^n \setminus F_{\lambda}$, we have 
	\[
		J_{h_\lambda} \leq \abs{Dh_{\R, \lambda}} \abs{Dh_{\S^{n-1}}}^{n-1} \leq C_n \lambda \abs{Dh}^{n-1}.
	\]
	Thus, we may estimate
	\begin{align*}
		\abs{\int_{F_\lambda} J_h}
		&= \abs{\int_{\R^n \setminus F_\lambda} J_{h_\lambda}}
		\leq C_n \int_{\R^n \setminus F_\lambda} \lambda \abs{Dh}^{n-1}. 
	\end{align*}
	Combining this with \eqref{eq:logdistineq}, we thus obtain
	\[
		\int_{F_\lambda} \abs{Dh}^n 
		\leq C_n K \int_{\R^n \setminus F_\lambda} \lambda \abs{Dh}^{n-1}
		+ \int_{F_\lambda} \sigma^n.
	\]
	Multiplying by $\Phi_{n}(\lambda)$, we get
	\begin{equation}\label{eq:li_step1}
		\int_{F_\lambda} \Phi_{n}(\lambda) \abs{Dh}^n 
		\leq C_n K \int_{\R^n \setminus F_\lambda} \lambda \Phi_{n}(\lambda) \abs{Dh}^{n-1}
		+ \int_{F_\lambda} \Phi_{n}(\lambda) \sigma^n .
	\end{equation}

	We then let $t > 0$, integrate \eqref{eq:li_step1} from $t$ to $\infty$ with respect to $\lambda$, and use Fubini-Tonelli to swap the order of integrals. Note that the set $F_\lambda$ is defined with the condition $M(\abs{Dh}) \leq \lambda$; in our use of Fubini-Tonelli, this condition must be moved to the integral with respect to $\lambda$, where it changes the lower bound of integration to $M(\abs{Dh})$ whenever $M(\abs{Dh}) > t$. Similarly, $\R^n \setminus F_\lambda$ is defined with the condition $M(\abs{Dh}) > \lambda$; when moved to the integral with respect to $\lambda$, this changes the upper bound to $M(\abs{Dh})$ when $M(\abs{Dh}) > t$, and makes the integral vanish otherwise. Altogether, we obtain
	\begin{multline}\label{eq:li_step2}
		\int_{\R^n} \abs{Dh}^n \left( \int_{\max(t, M(\abs{Dh}))}^\infty \Phi_{n}(\lambda) \dd\lambda \right)\\
		\leq C_n K \int_{\R^n}  \abs{Dh}^{n-1} \left(\int_{t}^{\max(t, M(\abs{Dh}))} \lambda \Phi_{n}(\lambda) \dd \lambda \right)\\
		+ \int_{\R^n} \sigma^n \left(\int_{\max(t, M(\abs{Dh}))}^\infty \Phi_{n}(\lambda) \dd\lambda \right).
	\end{multline}

	We then evaluate the integrals and estimate. For the $\sigma$-term, we use \eqref{eq:Phi_integral_identity} and the fact that we could assume that $\sigma \leq (1+K)^{1/n} \abs{Dh}$. Thus, in the region $\{\sigma \neq 0\}$, we get that
	\begin{multline*}
		\int_{\max(t, M(\abs{Dh}))}^\infty \Phi_{n}(\lambda) \dd\lambda 
		= \log^{n}\left( 1 + \frac{1}{\max(t, M(\abs{Dh}))} \right)\\
		\leq \log^{n}\left( 1 + \frac{1}{M(\abs{Dh})}\right)
		\leq \log^{n}\left( 1 + \frac{1}{\abs{Dh}}\right)
		\leq \log^{n}\left( 1 + \frac{(1+K)^{\frac{1}{n}}}{\sigma}\right).
	\end{multline*}
	Moreover, we recall that for all $a, b, \gamma > 0$, we have the elementary inequalities $\log(a) \leq C(\gamma) a^\gamma$ and $(a+b)^\gamma \leq C(\gamma)(a^\gamma + b^\gamma)$. By applying these inequalities in the region $\{\sigma \neq 0\}$, we obtain the further estimate
	\[
		\log^{n}\left( 1 + \frac{(1+K)^{\frac{1}{n}}}{\sigma}\right) \leq C(n, K, \eps) \left( 1 + \sigma^{-\eps} \right)
	\]
	In the region $\{\sigma = 0\}$, it suffices to know that since $M(\abs{Dh}) > 0$ everywhere due to $h$ being non-constant, the integral with respect to $\lambda$ in the $\sigma$-term is finite and is hence eliminated by the $\sigma^n$-coefficient. 
	
	For the $\abs{Dh}$-term on the left hand side, we similarly use \eqref{eq:Phi_integral_identity}, obtaining
	\begin{multline*}
		\int_{\max(t, M(\abs{Dh}))}^\infty \Phi_{n}(\lambda) \dd\lambda 
		= \log^{n}\left( 1 + \frac{1}{\max(t, M(\abs{Dh}))} \right)\\
		\geq \chi_{\R^n \setminus F_t} \log^{n}\left( 1 + \frac{1}{M(\abs{Dh})} \right).
	\end{multline*}
	The remaining middle integral is then estimated by using Lemma \ref{lem:Phi_integral_lemma}:
	\begin{multline*}
		\int_{t}^{\max(t, M(\abs{Dh}))} \lambda \Phi_{n}(\lambda) \dd \lambda 
		\leq \chi_{\R^n \setminus F_t} \int_{0}^{M(\abs{Dh})} \lambda \Phi_{n}(\lambda) \dd \lambda\\
		\leq \chi_{\R^n \setminus F_t} C(n) M(\abs{Dh}) \left( 1 + \log^{n-1}\left( 1 + \frac{1}{M(\abs{Dh})} \right) \right).
	\end{multline*}
	By applying all these estimates to \eqref{eq:li_step2}, we get
	\begin{multline}\label{eq:li_step3}
		\int_{\R^n \setminus F_t} \abs{Dh}^n \log^{n}\left( 1 + \frac{1}{M(\abs{Dh})} \right)\\
		\leq C(n) K\int_{\R^n \setminus F_t} \abs{Dh}^{n-1} M(\abs{Dh}) \log^{n-1}\left( 1 + \frac{1}{M(\abs{Dh})} \right) \\
		+ C(n) K \int_{\R^n \setminus F_t} \abs{Dh}^{n-1} M(\abs{Dh}) + C(n, K, \eps) \int_{\R^n} \left( \sigma^n + \sigma^{n-\eps} \right).
	\end{multline}
	
	We proceed to apply Young's inequality $\abs{ab} \leq \smallabs{a}^p/p + \smallabs{b}^{p^*}/p^*$ on the first term on the right hand side of \eqref{eq:li_step2}. In particular, we get
	\begin{multline*}
		C(n) K \abs{Dh}^{n-1} M(\abs{Dh}) \log^{n-1}\left( 1 + \frac{1}{M(\abs{Dh})} \right)\\
		\leq \frac{n-1}{n} \abs{Dh}^{n} \log^{n}\left( 1 + \frac{1}{M(\abs{Dh})} \right)
		+ \frac{(C(n)K)^n}{n} M^n(\abs{Dh}).
	\end{multline*}
	Thus, with a further use of the fact that $\abs{Dh} \leq M(\abs{Dh})$ a.e.\ in $\R^n$, we have
	\begin{multline}\label{eq:li_step4}
		\int_{\R^n \setminus F_t} \abs{Dh}^n \log^{n}\left( 1 + \frac{1}{M(\abs{Dh})} \right)\\
		\leq \frac{n-1}{n} \int_{\R^n \setminus F_t} \abs{Dh}^{n} \log^{n}\left( 1 + \frac{1}{M(\abs{Dh})} \right) \\
		+ C(n, K) \int_{\R^n \setminus F_t} M^n(\abs{Dh}) + C(n, K, \eps) \int_{\R^n} \left( \sigma^n + \sigma^{n-\eps} \right).
	\end{multline}
	We then claim that the first term on the right hand side of \eqref{eq:li_step4} is finite, and can hence be absorbed to the left hand side of \eqref{eq:li_step4}. Indeed, in $\R^n \setminus F_t$, we have $M(\abs{Dh}) > t$. Since $M(\abs{Dh}) \in L^n(\R^n)$ by the Hardy-Littlewood maximal inequality, the set $\R^n \setminus F_t$ must have finite measure. Hence, by again using the elementary estimate $\log(1+r) \leq r$ for $r \geq 0$, we have
	\[
		\int_{\R^n \setminus F_t} \abs{Dh}^n \log^{n}\left( 1 + \frac{1}{M(\abs{Dh})} \right)
		\leq \int_{\R^n \setminus F_t} \frac{\abs{Dh}^n}{M^n(\abs{Dh})} \leq \abs{\R^n \setminus F_t} < \infty.	
	\]
	Thus, subtracting the term from both sides is possible, leaving us with the estimate
	\begin{multline*}
		\frac{1}{n}\int_{\R^n \setminus F_t} \abs{Dh}^n \log^{n}\left( 1 + \frac{1}{M(\abs{Dh})} \right)\\
		\leq C(n, K) \int_{\R^n \setminus F_t} M^n(\abs{Dh}) + C(n, K, \eps) \int_{\R^n} \left( \sigma^n + \sigma^{n-\eps} \right).
	\end{multline*}
	
	It remains to estimate the $M^n(\abs{Dh})$-term using the global Hardy-Littlewood maximal inequality. In particular,
	\begin{multline*}
		\int_{\R^n \setminus F_t} M^n(\abs{Dh})
		\leq \int_{\R^n} M^n(\abs{Dh})
		\leq C(n) \int_{\R^n} \abs{Dh}^n\\
		\leq C(n) \int_{\R^n} (K J_h + \sigma^n)
		= C(n) \int_{\R^n} \sigma^n,
	\end{multline*}
	since the integral of $J_h$ over $\R^n$ vanishes by $\abs{Dh} \in L^n(\R^n)$ and \cite[Lemma 2.4]{Kangasniemi-Onninen_Heterogeneous} of the original article.  In conclusion,
	\[
		\int_{\R^n \setminus F_t} \abs{Dh}^n \log^{n}\left( 1 + \frac{1}{M(\abs{Dh})} \right)\\
		\leq C(n, K, \eps) \int_{\R^n} (\sigma^n + \sigma^{n-\eps}) < \infty.
	\]
	Notably, this bound is independent of $t$. Moreover, since $h$ is non-constant, we have $M(\abs{Dh}) > 0$ everywhere, and thus $\bigcup_{t > 0} \R^n \setminus F_t = \R^n$. We may thus let $t \to 0$ and use monotone convergence, and the claim follows.
\end{proof}

It now remains to prove Theorem 1.3. It turns out that the revised Lemma \ref{prop:log_integrability} is enough to complete a chain of balls estimate similar to the original proof. However, the basic sequence of balls with a shared center that was used in the original proof is barely not enough for a successful proof, and we instead need to use a sequence with limited overlaps that is closer in spirit to the chain condition of \cite{Hajlasz-Koskela-SobolevmetPoincare}.

\begin{customthm}{1.3}
	Suppose that $f \in W^{1,n}_\loc(\R^n, \R^n)$ satisfies the heterogeneous distortion inequality with $K \in [1, \infty)$ and $\sigma \in L^{n-\eps}(\R^n) \cap L^{n+\eps}(\R^n)$, for some $\eps > 0$. If $f$ is bounded and $\lim_{x \to \infty} \abs{f(x)} = 0$, then $f \equiv 0$.
\end{customthm}
\begin{proof}
	We start as in the original flawed proof. We suppose towards contradiction that $f$ is bounded and $\lim\limits_{x \to \infty} \abs{f(x)} = 0$, but $f$ is not identically zero. By \cite[Theorems 1.1 and 1.2]{Kangasniemi-Onninen_Heterogeneous} of the original article, $f$ is continuous and has no zeros. Thus, if we define $h \colon \R^n \to \R \times \S^{n-1}$ by
	\[
		h(x) = \left( \log \abs{f}, \frac{{f}}{\abs{f}} \right). 
	\]
	then \cite[Lemma 7.1]{Kangasniemi-Onninen_Heterogeneous} of the original article yields that $h \in W^{1,n}_\loc(\R^n, \R \times \S^{n-1})$, $\abs{Dh} \in L^n(\R^n)$, and $\abs{Dh}^n \leq K J_h + \sigma^n$ as before.
	
	We fix a unit vector $x_0 \in \S^{n-1} \subset \R^n$, and consider balls of the form $B_i = B^n(x_i, r_i)$, $i \in \Z_{> 0}$, where $x_i = (2^{i+1} - 1)x_0$ and $r_i = 2^{i}$. This selection ensures that $x_i$ is on the boundary of $B_{i+1}$, since $\abs{x_{i+1} - x_i} = (2^{i+2} - 1) - (2^{i+1} - 1) = 2^{i+1} = r_{i+1}$. It follows that for every $i$, the intersection of consecutive balls $B_i \cap B_{i+1}$ contains the ball $B'_i = B^n(x_i + (r_{i}/2) x_0, r_{i}/2)$. It is also crucial to observe that no point of $\R^n$ is contained in more than two of the balls $B_i$. See Figure \ref{fig:chain_of_balls} for an illustration of the balls $B_i$ and $B_i'$.
	
	\begin{figure}[h]
		\centering
		\begin{tikzpicture}[scale=0.3]
			\filldraw[color=gray, fill=gray!30] (4,0) circle (1);
			\filldraw[color=gray, fill=gray!30] (9,0) circle (2);
			\filldraw[color=gray, fill=gray!30] (19,0) circle (4);
			
			\draw (3,0) circle (2);
			\draw (7,0) circle (4);
			\draw (15,0) circle (8);
			\draw (15,0) arc (180:220:16);
			\draw (15,0) arc (180:140:16);
			
			\node at (3,-2)[anchor=north east] {$B_1$};
			\node at (7,-4)[anchor=north east] {$B_2$};
			\node at (15,-8)[anchor=north east] {$B_3$};
			
			\filldraw[black] (1,0) circle (0.2);
			\filldraw[black] (3,0) circle (0.2);
			\filldraw[black] (7,0) circle (0.2);
			\filldraw[black] (15,0) circle (0.2);
			
			\node at (1,0) [anchor=east] {$x_0$};
			\node at (3,0) [anchor=east] {$x_1$};
			\node at (7,0) [anchor=east] {$x_2$};
			\node at (15,0) [anchor=east] {$x_3$};
			
			\node at (4,0) [color=black!70, anchor=center] {$B_1'$};
			\node at (9,0) [color=black!70, anchor=center] {$B_2'$};
			\node at (19,0) [color=black!70, anchor=center] {$B_3'$};
		\end{tikzpicture}
		\caption{The points $x_i$ and the balls $B_i$ and $B_i'$.}\label{fig:chain_of_balls}
	\end{figure}
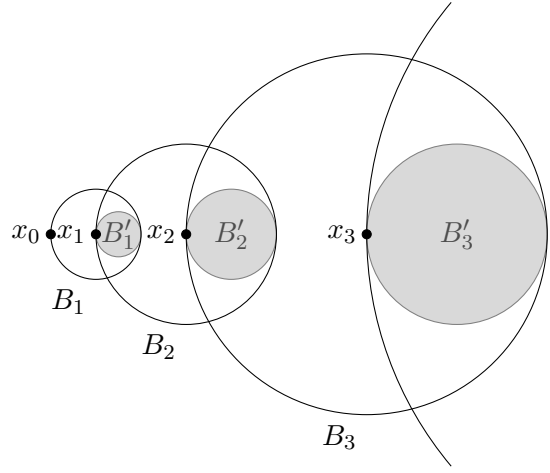
	
	Since $\log \abs{f} \in W^{1,n}_\loc(\R^n)$, it has a finite average integral $(\log \abs{f})_{B_{1}} \in \R$ over the first ball. Since $\lim_{x \to \infty} f(x) = 0$ and $\dist(0, B_i) \to \infty$, we also have $\lim_{i \to \infty} \max_{B_i} \log \abs{f} = -\infty$, and thus $\lim_{i \to \infty} (\log \abs{f})_{B_{i}} = -\infty$. Our objective is to show an $i$-independent upper bound for $\abs{(\log \abs{f})_{B_{i}} - (\log \abs{f})_{B_{1}}} = \abs{(h_\R)_{B_{i}} - (h_\R)_{B_{1}}}$, leading to a contradiction as in the original proof.
	
	We estimate this difference with a telescopic sum, obtaining
	\[
		\abs{(h_\R)_{B_{i}} - (h_\R)_{B_{1}}}
		\leq \sum_{j=1}^{i-1} \abs{(h_\R)_{B_{j+1}} - (h_\R)_{B_{j}}}
		\leq \sum_{j=1}^\infty \abs{(h_\R)_{B_{j+1}} - (h_\R)_{B_{j}}}.
	\]
	We then apply the Sobolev-Poincar\'e inequality, the fact that $B_j' \subset B_j \cap B_{j+1}$, and the fact that $r_{j+1}/r_{j} = 2$. As a result, we can estimate the terms in this telescopic sum by
	\begin{multline*}
		\abs{(h_\R)_{B_{j+1}} - (h_\R)_{B_{j}}}
		\leq \big\lvert(h_\R)_{B_j'} - (h_\R)_{B_{j}}\big\rvert
			+ \big\lvert(h_\R)_{B_{j}'} - (h_\R)_{B_{j+1}}\big\rvert\\		
		\leq \dashint_{B_{j}'} \abs{h_\R - (h_\R)_{B_{j}}} + \dashint_{B_{j}'} \abs{h_\R - (h_\R)_{B_{j+1}}}\\
		\leq 2^n \dashint_{B_{j}} \abs{h_\R - (h_\R)_{B_{j}}} + 4^n \dashint_{B_{j+1}} \abs{h_\R - (h_\R)_{B_{j+1}}}\\
		\leq C_n \left( r_j \dashint_{B_j} \abs{\nabla h_\R} + r_{j+1} \dashint_{B_{j+1}} \abs{\nabla h_\R} \right) .
	\end{multline*}
	In particular, since $\abs{\nabla h_\R} \leq \abs{Dh}$, we have the estimate
	\begin{equation}\label{eq:telescopic_term_estimate}
		\abs{(\log \abs{f})_{B_{i}} - (\log \abs{f})_{B_{1}}}
		\leq 2C_n V_n^{-1} \sum_{j=1}^\infty r_j^{-n-1} \int_{B_j} \abs{Dh},
	\end{equation}
	where $V_n$ denotes the volume of the unit ball in $\R^n$.
	
	We then recall the following elementary inequality: if $I \subset \R$ is a possibly infinite interval, $F, G \colon I \to (0, \infty)$ with $F$ increasing and $G$ decreasing, and $s, t \in I$, then
	\[
		1 \leq \max \left( \frac{F(t)}{F(s)}, \frac{G(t)}{G(s)} \right) \leq \frac{F(t)}{F(s)} + \frac{G(t)}{G(s)}.
	\]
	We apply this to the right hand side of \eqref{eq:telescopic_term_estimate} with $I = (0, \infty)$, $F(\rho) = \rho^{n-1}$, $G(\rho) = \log(1 + \rho^{-1})$, $t = M(\abs{Dh})$, and $s = r_j^{-1/2}$. In particular, we have
	\begin{multline}\label{eq:two_terms}
		\sum_{j=1}^\infty r_j^{-n-1} \int_{B_j} \abs{Dh}
		\leq \sum_{j=1}^\infty r_j^{-\frac{n-1}{2}} \int_{B_j} \abs{Dh} M^{n-1}(\abs{Dh})\\
		+ \sum_{j=1}^\infty \frac{r_j^{-n-1}}{\log(1 + \sqrt{r_j})} \int_{B_j} \abs{Dh} \log \left(1 + \frac{1}{M(\abs{Dh})} \right).
	\end{multline}
	For the first term on the right hand side of \eqref{eq:two_terms}, since $r_j = 2^j$ and since $\abs{Dh} \leq M(\abs{Dh})$ a.e.\ in $\R^n$, we may estimate that
	\[
		\sum_{j=1}^\infty r_j^{-\frac{n-1}{2}} \int_{B_j} \abs{Dh} M^{n-1}(\abs{Dh})
		\leq \left( \int_{\R^n} M^{n}(\abs{Dh}) \right) \sum_{j=1}^\infty 2^{-\frac{(n-1)j}{2}}.
	\]
	Since $\abs{Dh} \in L^n(\R^n)$, the integral of $M^n(\abs{Dh})$ over $\R^n$ is finite, and moreover the sum part of the above upper bound is a convergent geometric sum. Hence, we've obtained an upper bound independent of $i$ for this part of the sum.
	
	For the second term on the right hand side of \eqref{eq:two_terms}, we first apply H\"older's inequality for integrals with 1 as the other function, obtaining
	\begin{multline*}
		\sum_{j=1}^\infty \frac{r_j^{-n-1}}{\log(1 + \sqrt{r_j})} \int_{B_j} \abs{Dh} \log \left(1 + \frac{1}{M(\abs{Dh})} \right)\\
		\leq C(n) \sum_{j=1}^\infty \frac{1}{\log(1 + \sqrt{r_j})} \left( \int_{B_j} \abs{Dh}^n \log^n \left(1 + \frac{1}{M(\abs{Dh})} \right) \right)^\frac{1}{n}.
	\end{multline*}
	Next, we use H\"older's inequality for sums, and apply the fact that every point of $\R^n$ is an element of at most two of the sets $B_j$. The resulting estimate is
	\begin{multline*}
		\sum_{j=1}^\infty \frac{1}{\log(1 + \sqrt{r_j})} \left( \int_{B_j} \abs{Dh}^n \log^n \left(1 + \frac{1}{M(\abs{Dh})} \right) \right)^\frac{1}{n}\\
		\leq \left( \sum_{j=1}^\infty \frac{1}{\log^\frac{n}{n-1}(1 + \sqrt{r_j})} \right)^\frac{n-1}{n}
		\left( \sum_{j=1}^\infty \int_{B_j} \abs{Dh}^n \log^n \left( 1 + \frac{1}{M(\abs{Dh})}\right) \right)^\frac{1}{n}\\
		\leq \left( \sum_{j=1}^\infty \frac{1}{\log^\frac{n}{n-1}(1 + \sqrt{r_j})} \right)^\frac{n-1}{n}
		\left( 2 \int_{\R^n} \abs{Dh}^n \log^n \left( 1 + \frac{1}{M(\abs{Dh})}\right) \right)^\frac{1}{n}.
	\end{multline*}
	Now, by the revised Lemma \ref{prop:log_integrability}, the integral term in this estimate is finite. On the other hand, for the sum term, we get by $r_j = 2^j$ that
	\[
		\sum_{j=1}^\infty \frac{1}{\log^\frac{n}{n-1}(1 + \sqrt{r_j})}
		< \sum_{j=1}^\infty \frac{1}{\log^\frac{n}{n-1}(\sqrt{r_j})}
		= \left(\frac{2}{\log(2)}\right)^\frac{n}{n-1} \sum_{j=1}^\infty j^{-\frac{n}{n-1}} < \infty
	\]
	since $n/(n-1) > 1$. Thus, we have bounded the right hand side of \eqref{eq:telescopic_term_estimate} with a finite bound independent of $i$, which is a contradiction. The claimed result follows.
\end{proof}

\bibliographystyle{abbrv}
\bibliography{sources}

\end{document}